\scrollmode 
\documentclass[12pt]{amsart} 
\usepackage{verbatim}
\usepackage{eucal} 
\usepackage{psfrag}
\usepackage{graphicx}
\usepackage{mathrsfs}
\theoremstyle{plain}
\newtheorem{theorem}{Theorem}
\newtheorem{proposition}{Proposition}
\newtheorem{lemma}{Lemma}

\theoremstyle{definition}
\newtheorem{definition}{Definition}
\newtheorem{example}{Example}

\theoremstyle{remark}
\newtheorem{remark}{Remark}

\newcommand{\R}{{\mathbb R}}

\newcommand{\C}{{\mathbb C}}
\newcommand{\Z}{{\mathbb Z}}
\newcommand{\Q}{{\mathbb Q }}

\newcommand{\G}{{\mathscr G}}
\newcommand{\OK}{{\mathscr{O}_{K}}}

\newcommand{\hp}{{\mathfrak H}}
\newcommand{\Vor}{Vorono\v{\i}}
\newcommand{\Int}{\operatorname{Int}}

\newcommand{\Hom}{\operatorname{Hom}}
\newcount\timehh\newcount\timemm
\timehh=\time 
\divide\timehh by 60 \timemm=\time
\newcommand{\draftauthor}[1]{\author{#1
    {
      --- \protect \protect\sc\today\ ---
      \ifnum\timehh<10 0\fi\number\timehh\,:\,\ifnum\timemm<10 0\fi\number\timemm
      \protect \, \, \protect \bf DRAFT
    }
  }
}

\begin{document}
%
%
\title{Modular Symbols for $\Q $-rank one groups and Vorono\v{\i}\  Reduction}
\author{Paul E. Gunnells}
\address{Department of Mathematics\\
Columbia University\\
New York, New York  10027}
\date{\today}
\email{gunnells@math.columbia.edu}
%
%
\begin{abstract}
Let $G$ be a reductive algebraic group of $\Q $-rank one associated to
a self-adjoint homogeneous cone defined over $\Q $, and let $\Gamma
\subset G$ be a torsion-free arithmetic subgroup.  Let $d$ be the
cohomological dimension of $\Gamma $.  We present an algorithm to
compute the action of the Hecke operators on $H^{d}(\Gamma ;\Z )$.
This generalizes the classical modular symbol algorithm, when $\Gamma
\subset SL_{2}(\Z )$, to a setting including Bianchi groups and
Hilbert modular groups.  In addition, we generalize some results of
\Vor\ for real positive-definite quadratic forms to 
self-adjoint homogeneous cones of arbitrary $\Q $-rank.
\end{abstract}
\maketitle
%
%
\section{Introduction}\label{introduction}
\subsection{}\label{background} 
Let $\hp _{2} = SL_{2}(\R )/SO(2)$ be the upper-half plane, and let
$\Gamma(N) \subset SL_{2}(\Z )$ be the principal congruence subgroup
of level $N>2$, acting on $\hp _{2}$ from the left by
linear fractional transformations.  Then the cohomology group
$H^{1}(\Gamma(N)\backslash \hp _{2};\C )$ is closely related to the
space of all weight-two modular forms of level $N$.  The {\em modular
symbols} provide a concrete approach to the group
$H^{1}(\Gamma(N)\backslash \hp _{2};\C )$ (\S\ref{manin.work}) that
has allowed the testing of many conjectures in number theory and has
led to explicit formulas for $L$-functions and their derivatives
\cite{diamantis}\cite{goldfeld}\cite{manin}.  Important to
applications is the modular symbol algorithm developed by Manin
\cite{manin}.  An algebra of Hecke operators acts on
$H^{1}(\Gamma(N)\backslash \hp _{2};\C )$, and using the algorithm one
may compute their eigenvalues.  Essentially this algorithm is the
euclidean algorithm applied to pairs of integers
(\S\ref{classicalmodularsymb}). 

Now consider the case where $\Gamma $ is a torsion-free Bianchi
subgroup, that is, $\Gamma$ is of finite index in $SL_{2}(\OK)$, where
$\OK$ is the ring of integers in an imaginary quadratic extension
$K/\Q $.  The group $\Gamma $ acts on hyperbolic three-space $\hp _{3}
= SL_{2}(\C)/SU(2)$, and we consider the cohomology group
$H^{2}(\Gamma \backslash \hp _{3};\C )$.  As before, an algebra of
Hecke operators acts on the cohomology, and one is interested in this
Hecke-module for many reasons.  For example, results of Grunewald and
Schwermer~\cite{grune.sch} imply that, for all but a finite set of
$K$, the rational cohomology of $SL_{2}(\OK)$ contains cuspidal
cohomology, which is important in the theory of automorphic forms.
Also, the ``Langlands philosophy'' predicts a
``Shimura-Taniyama-Weil'' correspondence between Hecke eigenclasses
and certain algebraic varieties defined over $K$.  More precisely, let
$\Gamma \subset SL_{2}(\OK)$ be a congruence subgroup and let $\xi \in
H^{2}(\Gamma \backslash \hp _{3};\C )$ be a cuspidal Hecke eigenclass.
Then one hopes to associate to $\xi $ an algebraic variety
$V/K$---specifically an elliptic curve or an abelian variety of
dimension two---so that the zeta function of $V$ is assembled from the
eigenvalues of $\xi $ in a precise way.  Results of Cremona
\cite{crem} and Cremona and Whitley~\cite{crem.whit} when $K$ has
class number one support this.  Hence one wishes to compute Hecke
eigenvalues for general $K$.  But in general the ring $\OK$ is not a
euclidean domain, and one cannot directly apply the modular
symbol algorithm as described in~\cite{manin} (however, see
\S\ref{related}).

In this paper we present an analog of the modular symbol algorithm for
$H^{d}(\Gamma;\Z )$, where $\Gamma $ is a torsion-free arithmetic
group associated to certain self-adjoint homogeneous cones, and $d$ is
the cohomological dimension of $\Gamma $.  This includes finite index
subgroups of $SL_{2}(R)$, where $R$ is
\begin{samepage}
\begin{itemize}
\item $\Z $,
\item the ring of integers in a $CM$ field, or 
\item the ring of integers in a totally real field.
\end{itemize}
\end{samepage}
We replace the continued fractions of~\cite{manin} with a study of the
geometry of self-adjoint homogeneous cones.  Thus our algorithm does
not require that $R$ be a euclidean domain.

Here is the organization of this paper.  Section~\ref{example}
contains a review of the classical modular symbol algorithm from
\cite{manin} and presents our algorithm in that case.  Section
\ref{red.thy} contains a review of the reduction theory of
self-adjoint homogeneous cones, and in Section~\ref{vor.red} we
generalize some results of~\cite{voronoi2} to this setting (Theorems
\ref{finite.num.of.steps} and~\ref{reduction.alg}).  These results are
valid for cones of any $\Q $-rank $\geq 1$.  Finally, Section
\ref{mod.sym} contains a description of our algorithm, in Theorem
\ref{hecke.alg}.  Throughout the paper we comment on implementation
issues related to the algorithm.

\subsection{Related work}\label{related}
Let $\Gamma \subset SL_{2}(\OK)$, where $\OK$ is the ring of
integers in an imaginary quadratic extension $K = \Q (\sqrt{-m})$.
For $K$ a non-euclidean ring with class number one
($m=19,43,67,163$), Whitley developed a ``pseudo-euclidean algorithm''
that allowed implementation of the modular symbol algorithm
\cite{whitley}.

Also, I learned upon completion of this work that Jeremy Bygott
has independently studied the modular symbol algorithm
for the non-PID imaginary quadratic case in his forthcoming
Ph.D. thesis~\cite{bygott}.

\subsection{Acknowledgments}

The results in this paper depend heavily on results from the work of
Avner Ash.  I thank him for graciously and patiently
explaining his work to me.  I also thank Mark McConnell
for a careful reading of an early version of this paper and many
helpful comments.  I also thank the referee for many
helpful suggestions.

Finally, the results in this paper are an extension of some of
the results in my Ph.D. thesis~\cite{me}.  I thank heartily my
advisor, Robert MacPherson, for the encouragement and inspiration he
has given me.

\section{A motivating example}\label{example}
In this section we illustrate our algorithm for $\Gamma
\subset SL_{2}(\Z )$.  As in the introduction let $\hp _{2} =
SL_{2}(\R )/SO_{2}$, and let $\hp _{2}^{*} = \hp _{2}\cup \Q
\cup \{\infty \}$ be the usual partial compactification of $\hp
_{2}$ by adding cusps.  We assume that $\hp _{2}^{*}$ is given the
Satake topology, and we extend the action of $SL_{2}(\R )$ to the
cusps.  We denote the quotient $\Gamma \backslash \hp _{2}^{*}$ by
$X_{\Gamma }$.  We assume that $\Gamma $ is torsion-free, so that
$X_{\Gamma }$ is smooth.

\subsection{}\label{manin.work}
We begin by paraphrasing aspects of Manin's work~\cite{manin}.  By
Poincar\'e duality, $H^{1}(X_{\Gamma} ;\C ) $ may be identified with
$H_{1}(X_{\Gamma};\C )$, and so we may study the space of weight-two
modular forms for $\Gamma $ by studying the latter.  Let
$q_{1}$ and $q_{2}$ be cusps equivalent modulo $\Gamma $.  Then
any smooth path $\gamma $ from $q_{1}$ to $q_{2}$ descends to a closed
path on $X_{\Gamma }$ representing a class in $H_{1}(X_{\Gamma };\Z
)$.  Furthermore, this class is independent of $\gamma $, and in fact
depends only on the ordered pair $(q_{1},q_{2})$.

More generally, suppose $q_{1}$ is not necessarily equivalent to
$q_{2}$ modulo $\Gamma $.  Then integration of one-forms $\omega \in
H^{1}(X_{\Gamma };\R )$ along $\gamma $ yields a functional $\int
\colon H^{1}(X_{\Gamma };\R )\rightarrow \R $, and this allows us
to associate to the pair $(q_{1},q_{2})$ a class in $H_{1}(X_{\Gamma
};\R )$.  By the theorem of Manin-Drinfeld (\cite{lang}, p.~61), this
class actually lies in $H_{1}(X_{\Gamma };\Q )$.  We define a
{\em modular symbol} to be the rational homology class constructed
from an ordered pair of cusps in this way, and denote this class by
$[q_{1},q_{2}]$.  This class agrees with the class in the previous
paragraph when $q_{1}$ and $q_{2}$ are equivalent modulo $\Gamma $.

\begin{proposition}\label{manin.sl2}
\cite{manin} The modular symbols satisfy the following:
\begin{enumerate}
\item $[q_{1},q_{2}] = - [q_{2},q_{1}]$.
\item $[q_{1},q_{2}] = [q_{1},q_{3}] + [q_{3},q_{2}]$. 
\end{enumerate}
Furthermore, $H_{1}(X_{\Gamma };\Z )$ is spanned by modular symbols
modulo $\Gamma $.
\end{proposition}

\subsection{}\label{classicalmodularsymb}
The Hecke operators act on $H^{1}(X_{\Gamma };\C) $, and hence by
duality on $H_{1}(X_{\Gamma };\C) $.  On a modular symbol,
an operator acts by 
$$[q_{1},q_{2}]\longmapsto \sum _{\alpha \in A} \,[\alpha q_{1},\alpha
q_{2}], $$ 
where $A$ is a finite set of $2\times 2$
integral matrices that depends on the operator.  For example, let 
$\Gamma =\Gamma (N)$, and let $p$ be a prime not dividing $N$.  For
the 
classical operator $T_{p}$ we may take
$$\alpha \in \left\{%
	\sigma _{p}\left( \begin{array}{cc}
         p&0\\
         0&1	  
	\end{array}\right),
	\left( \begin{array}{cc}
         1&0\\
         0&p	  
	\end{array}\right),
	\left( \begin{array}{cc}
         1&1\\
         0&p	  
	\end{array}\right),
        \ldots,
	\left( \begin{array}{cc}
         1&p-1\\
         0&p	  
	\end{array}\right)
 \right\}.%
 $$
Here $\sigma _{p}\in SL_{2}(\Z )$ is a fixed matrix satisfying 
$$\sigma _{p}\equiv \left(\begin{array}{cc}
p^{-1}&0\\
0&p
\end{array} \right)\mod N $$
(\cite{shimura}, Prop.~3.36).  Note that $\det \alpha \not = \pm 1$.

A finite basis of $H_{1}(X_{\Gamma };\C)$ 
is provided by the set of {\em unimodular symbols}.  Write
a cusp 
$q$ in lowest terms as $m/n$, where the cusp $\infty $
is written formally as $1/0$.  Then the unimodular
symbols are the symbols $[q_{1},q_{2}]$ satisfying 
$$\det \left(\begin{array}{cc}
m_{1}&m_{2}\\
n_{1}&n_{2}
\end{array}\right) = \pm 1. $$

\begin{figure}[h]
\begin{center}
\psfrag{0}{$0$}
\psfrag{12/5}{$12/5$}
\psfrag{2}{$2$}
\psfrag{5/2}{$5/2$}
\includegraphics{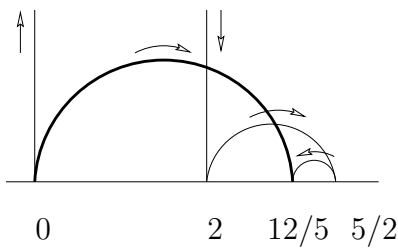}
\end{center}
\caption{\label{nonunimod} $12/5 = [\![2,2,2]\!]$ implies $[0,12/5] =
[0,\infty ] + [\infty ,2] + [2,5/2] + [5/2,12/5]$}
\end{figure}

The Hecke operators do not preserve unimodularity, and it is necessary
for eigenvalue computations to construct an explicit homology between
a non-unimodular symbol and a cycle of unimodular symbols.  This is
done by the {\em modular symbol algorithm}.  Assume
that a non-unimodular symbol has the form $[0,q]$, where $q$ is a
positive rational number.  Let $[\![a_{1},\ldots,a_{k}]\!]$ be the
simple continued fraction expansion of $q$, i.e.
$$q = a_{1} + \cfrac{1}{a_{2} + \cfrac{1}{\cdots +\cfrac{1}{a_{k}}}} $$
Let $q_{i}$ be the $i^{{\rm th}}$ convergent $[\![a_{1},\ldots,a_{i}]\!]$.
Then by applying (2) from Proposition~\ref{manin.sl2}, we have 
$$[0,q] = [0,\infty ] + [\infty ,q_{1}] + \cdots + [q_{k-1},q]. $$
Furthermore, the basic properties of simple continued fractions imply that 
the modular symbols on the right are unimodular.  Figure~\ref{nonunimod}
illustrates the result for the modular symbol $[0,12/5]$.

To complete the discussion, we note that $SL_{2}(\Z )$ acts
transitively on the cusps.  Since the above algorithm is $SL_{2}(\Z
)$-equivariant, any modular symbol can be written as a
sum of unimodular symbols.  

\subsection{}\label{example.myalg}
Now we present our technique for writing a modular symbol as an
equivalent sum of unimodular symbols.  No use will be made of
continued fractions; instead, we look at the relationship between a
geodesic representing a modular symbol and a certain tessellation of
$\hp _{2}$.  In this simple case our algorithm will appear needlessly
complicated, but it is formulated in a way that will generalize to
other settings.  It is also quite practical for machine computations.

\begin{figure}[h]
\begin{center}
\psfrag{inf}{$\infty $}
\psfrag{1}{$1$}
\psfrag{0}{$0$}
\psfrag{2}{$2$}
\includegraphics{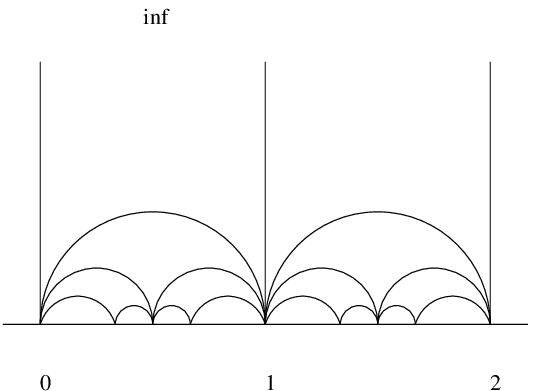}
\end{center}
\caption{\label{tess} Ideal triangles in $\hp _{2}$}
\end{figure}

To begin, we tile $\hp_{2}$ with the $SL_{2}(\Z )$-translates of the
ideal geodesic triangle with vertices $0$, $1$, and $\infty $ (see
Figure~\ref{tess}).  This tessellation descends to a finite
triangulation of $X_{\Gamma }$.  The edges of this tessellation are
geodesics inducing the unimodular symbols, and every unimodular symbol
arises in this way.

Given any point $x\in \hp _{2}$, let $R(x)$ be the set of
vertices of the triangle (or edge) of the tessellation meeting $x$.

\begin{figure}[h]
\centerline{ \includegraphics{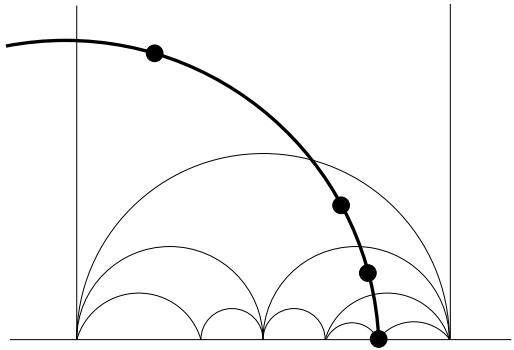} }
\caption{\label{partition} A partition of $\gamma $}
\end{figure}

Let $[0,q]$ be a modular symbol as before, and let $\gamma $ be
the ideal geodesic in $\hp _{2}$ from $0$ to $q$. 
Because $\gamma $ is a geodesic between two rational cusps, one can
show that $\gamma $ will
only meet a finite number of triangles in the tessellation.
Hence we may
choose a finite subset $x_{1},\ldots,x_{r}\in \gamma $ (as in Figure
\ref{partition}) so that
\begin{enumerate}
\item $0\in R(x_{1})$,
\item $q\in R(x_{r})$, and 
\item $R(x_{i})\cap R(x_{i+1})\not = \emptyset $.
\end{enumerate}

\noindent We call such a collection a {\em sufficiently fine partition} of
$\gamma $.  From each $R(x_{i})\cap R(x_{i+1})$ choose a cusp
$q_{i}$.  Then we claim that we have a homology
\begin{equation}\label{hom}
[0,q]=[0,q_{1}]+[q_{1},q_{2}]+\cdots +[q_{r},q],
\end{equation}
and that each term on the right is a unimodular symbol.

\begin{figure}[h]
\centerline{ \includegraphics{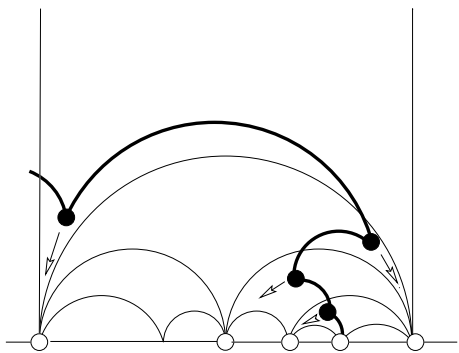} }
\caption{\label{deformation} A deformation of $\gamma $}
\end{figure}

First, we may see we have a homology either by repeatedly
applying (2) of Proposition~\ref{manin.sl2}, or by continuously
deforming $\gamma $ into geodesics inducing the classes on the
right-hand side (see Figure~\ref{deformation}).

Finally, each nontrivial term on the right of \eqref{hom} corresponds
to an edge in the tessellation because $q_{i}$ and $q_{i+1}$ are both
vertices of a triangle containing $x_{i+1}$.  In fact, as shown in
\cite{me}, if $\gamma $ is oriented then one may choose the $q_{i}$
canonically, and our algorithm in this case is equivalent to the
modular symbol algorithm.

\section{Self-adjoint homogeneous cones}\label{red.thy}
In this section we present the geometric context of our algorithm.  This specifies the arithmetic groups
to which our algorithm applies and describes the constructions
replacing $\hp _{2}$ and its tessellation from \S\ref{example}.  

The results in \S\ref{bkbegin}--\ref{bkend} and \S\ref{bkpoly} are due to
A. Ash and originally appeared in~\cite{ash} and~\cite{amrt}.  Our 
exposition closely follows the former.  
 
\subsection{} \label{bkbegin}
Let $V$ be a real vector space defined over $\Q $, and
let $C\subset V$ be an open cone.  That is, $C$ contains no straight
line, and $C$ is closed under homotheties: if $x\in C$ and $\lambda
\in \R^{>0}$, then $\lambda x\in C$.  The cone $C$ is
called {\em self-adjoint} if there exists a scalar product $\langle\,
,\, \rangle$ on $C$ such that
$$C = \left\{\,x\in V \,|\, \langle y,x\rangle > 0 \quad \hbox{for all $y\in \bar C - \{0 \}$} \,\right\}. $$

Let $G$ denote the connected component of the identity of the linear
automorphism group of $C$, i.e. $G=\left\{\,g\in GL(V)\,|\,gC = C\,
\right\}^{0}$.  The cone $C$ is called {\em homogeneous} if $G$ acts
transitively on $C$.  If $K$ denotes the isotropy group of a fixed
point in $C$, then we may identify $C$ with $G/K$.  The fact that $C$
is self-adjoint implies that $G$ is reductive and $C$ modulo
homotheties is a Riemannian symmetric space.

We also assume that all these notions are compatible with the $\Q
$-structure on $V$.  That is, as a subgroup of $GL(V)$, $G$ is defined
by rational equations, and the scalar product $\langle\, ,\, \rangle$
is defined over $\Q $.  This is stronger than saying that $G$ is
defined over $\Q $.  In particular, the group of real points $G (\R )$
must be isomorphic to a product of the following groups (\cite{cones},
p. 97):
\begin{enumerate}
\item $GL_{n} (\R )$
\item $GL_{n} (\C )$
\item $GL_{n} (\mathbb H )$
\item $O(1,n-1)\times \R ^{\times }$
\item The noncompact Lie group with Lie algebra $\mathfrak e_{6
(-26)}\oplus \R $ 
\end{enumerate}
In each case $V$ is a set of hermitian symmetric matrices.  In other
words, $V$ is the set of $n\times n$ matrices over an appropriate $\R
$-algebra with involution $\tau $, in which $A\in V$ if and only if
$A^{t}= A^{\tau }$.  The cone $C$ is then the subset of
``positive-definite'' matrices in an appropriate sense.  For details
we refer to (\cite{cones}, Ch. V).

\subsection{}
Let $H$ be a hyperplane in $V$.  We say that $H$ is a {\em supporting
hyperplane} of $C$ if $H$ is rational and $H\cap C=\emptyset $ but
$H\cap \bar C \not =\emptyset $.

Given a supporting hyperplane $H$ of $C$, let $C' =
\Int(H\cap \bar C)$.  (Here $\Int(A )$ is the interior of $A$
in its linear span.)  Then $C'$ is called
a {\em rational boundary component}, and is a
self-adjoint homogeneous cone of smaller dimension than $C$.

\begin{definition}\label{cusps}
The {\em cusps} of $C$ are the one-dimensional rational boundary
components of $C$.  The set of cusps is denoted $\Xi (C)$.
\end{definition}

\subsection{}
Let $L\subset V(\Q )$ be a lattice, i.e. a discrete subgroup of $V(\Q
)$ such that $L\otimes \Q =V(\Q )$.  Let $\Gamma_{L}$ denote the
subgroup of $G(\Q )$ preserving $L$.  An {\em arithmetic subgroup } of $G$
is a discrete subgroup commensurable with $\Gamma_{L}$ for some $L$.
Any torsion-free subgroup $\Gamma \subset \Gamma_{L}$ of finite index
will act properly discontinuously and freely on $C$.  Thus the
quotient $\Gamma \backslash C$ is an Eilenberg-Mac~Lane space for
$\Gamma $, and the group cohomology $H^{*}(\Gamma )$ is $
H^{*}(\Gamma\backslash C )$.  In fact, since homotheties commute with
the action of $\Gamma $, we may pass to $X := \R^{>0}\backslash C$,
and compute $H^{*}(\Gamma \backslash X)$ instead.

\subsection{}\label{bkend}
Let $A\subset V (\Q )$ be a finite set of nonzero points.  The closed
convex hull $\sigma $ of the rays $\left\{\R ^{\geq 0}x \bigm| x\in A
\right\}$ is called a {\em rational polyhedral cone}.  The rays
through the vertices of the convex hull of $A$ are called the {\em
spanning rays} of $\sigma $.  We denote the set of spanning rays by
$R(\sigma )$.  The group $G (\Q )$ acts naturally on the set of
rational polyhedral cones, and we denote the action by a dot: $\sigma
\mapsto g\cdot \sigma $.

Now we would like to partition $C$ into convex subsets using a
collection of rational polyhedral cones, in a manner compatible with
the $\Gamma_{L}$-action.  This requires some care, as $C$ is open
and any $\sigma $ as above is closed.

\begin{definition}(\cite{amrt}, p. 117)\label{gam.admiss}
Let $\Gamma \subset \Gamma _{L}$ be an arithmetic subgroup of $G$.
A set of closed polyhedral cones $\{\sigma _{\alpha } \}$ is called a
{\em $\Gamma $-admissible decomposition} of $C$ when the following hold:
\begin{enumerate}
\item Each $\sigma _{\alpha }$ is the span of a finite number of
rational rays.
\item For each $\alpha $, the cone $\sigma _{\alpha }\subset \bar C$.
\item Every face of a $\sigma _{\alpha }$ is some $\sigma _{\beta }$
in the decomposition.
\item $\sigma _{\alpha }\cap \sigma _{\beta }$ is a common face of
$\sigma _{\alpha } $ and $\sigma _{\beta }$.
\item For any $\sigma _{\alpha }$ and any $\gamma \in \Gamma $,
$\gamma \sigma _{\alpha }$ is some $\sigma _{\beta }$ in the decomposition.
\item Modulo $\Gamma$, there are only a finite number of $\sigma _{\alpha }$'s.
\item $C = \bigcup _{\alpha } (\sigma _{\alpha }\cap C)$.
\end{enumerate}
\end{definition}
Note that a $\Gamma $-admissible decomposition descends to a
decomposition of $X$ into open cells.

We now describe a technique to construct $\Gamma $-admissible
decompositions.  The technique originates with \Vor , and was
generalized by Ash to all self-adjoint homogeneous cones.  Let $L'$ be
$L-\{0 \}$.

\begin{definition}\label{vor.poly}
The {\em \Vor \ polyhedron} $\Pi $ is the closed convex hull of
$L'\cap \,\Xi (C)$.
\end{definition}

\begin{theorem}
(\cite{amrt}, p. 143) The cones over faces of $\Pi $ form a $\Gamma
$-admissible decomposition of $C$.
\end{theorem}

We call the cones in this $\Gamma $-admissible decomposition of $C$ the {\em
\Vor \ cones}, the decomposition of $X$ associated to
$\Pi $ the {\em \Vor \ decomposition}, and the open cells in $X$
the {\em \Vor \ cells}.\footnote{Our terminology is
nonstandard.  In~\cite{amrt}, the polyhedron $\Pi $ is referred to as
a ``kernel comparable to a $\Gamma $-polyhedral cocore'' and is a specific
example of a more general theory.}  Here are two examples of this
construction.

\begin{example}\label{classical}
The original example investigated by \Vor \ ~\cite{voronoi2} is the
following.  Let $V$ be the vector space of symmetric $n\times n$
matrices, and let $C\subset V$ be the cone of positive-definite
matrices.  Then $G$ is $GL_{n}(\R )^{0}$, which acts on $V$ by
$A\mapsto gAg^{t}$, and $K$ is $SO_{n}$.  The scalar product is given
by $\langle A, B \rangle = \hbox{Tr}(AB)$.

Let $L$ be the lattice of integral symmetric matrices.
Then $\Gamma_{L} = SL_{n}(\Z )$.  The set of cusps $\Xi (C)$ is
obtained as follows.  Any nonzero integral (column) $n$-vector $v$
determines a rank one quadratic form by $v\mapsto vv^{t}$.  The
cusps are the rays generated by all points in $\bar C$ of this form.
Suppose that $z\in L'\cap \Xi (C)$ is a cusp arising in this manner
from the $n$-vector $v$.  Then if $y\in C$, the scalar product
$\langle z, y \rangle$ is equal to the quadratic form $y$ evaluated on
$v$.

For $n=2$ we have that $X=\hp _{2}$, and the \Vor \ decomposition is
the tessellation described in \S\ref{example.myalg}.  For $n\geq 4$
there is more than one type of top-dimensional \Vor \ cone modulo
$\Gamma_{L} $, and not all the top-dimensional cones are
simplicial.  Complete information about these decompositions for $n=3$
and $4$ can be found in~\cite{mcc-macph1} and~\cite{mcc.art}.
\end{example}

\begin{example}\label{imag.quad} 
Let $K/\Q $ be an imaginary quadratic extension, and let $\OK$ be the
ring of integers of $K$.  Let $V$ be the vector space of $2\times 2$
hermitian symmetric matrices over $\C $, and let $C\subset V$ be the
cone of positive-definite matrices.  Then $V$ is a four-dimensional
vector space over $\R $, and $X$ is three-dimensional hyperbolic space
$\hp _{3}$.  For $L$ we may take the matrices in $V$ with entries in
$\OK$, and then $\Gamma_{L}=SL_{2}(\OK)$.  In the classical picture of
$\hp _{3}\subset \C \times \R^{\geq 0}$, the rays generated by the
vertices of $\Pi $ become the points $K\cup \{\infty \}$, where $K$
is pictured as a subset of $\C\times \{0 \} $ and $\infty $ is
pictured infinitely far above $\C \times \{0 \}$ along $\R^{\geq 0}$.

The \Vor \ decomposition becomes a tessellation of $\hp _{3}$ into
ideal three-polytopes.  In general these polytopes will not be
simplices.  For example, if $K=\Q (\sqrt{-1})$, then the unique
polytope modulo $\Gamma_{L} $ is an octahedron
\cite{crem}\cite{grune.hell.menn}.  If $\OK$ is euclidean there is
only one type of top-dimensional polytope in the tessellation, but for
general $K$ there will be more than one type.
\end{example}

\subsection{}\label{connection}
Here is the connection between the \Vor \ decomposition of $X$ and
$H^{*}(\Gamma ;\Z )$.  Let $N$ be the dimension of $X$ and $d$ the
cohomological dimension of $\Gamma $.  Let $C^{k}$ be the set of \Vor
\ cells of codimension~$k$.  The group $\Gamma $ acts naturally on
$C^{k}$ by its action on rational polyhedral cones.  We want to
construct a $\Gamma _{L}$-equivariant ``coboundary'' map $\delta
^{k}\colon \Z [C^{k}]\rightarrow \Z [C^{k+1}]$ so that the resulting
chain complex modulo $\Gamma $ computes $H^{*}(\Gamma;\Z )$.  We will
call $(C^{*},\delta ^{*})$ a {\em cocell complex}, and will say that
the \Vor \ decomposition gives $X$ a {\em cocell structure}.

According to~\cite{ash}, there is a topological space $W\subset X$
such that the following hold:
\begin{enumerate}
\item $W$ admits the structure of a $\Gamma $-equivariant regular cell
complex with top-dimensional cells of dimension $d$.
\item $W$ is a deformation retract of $X$, so that the homology of the
 chain complex associated to $\Gamma \backslash W$ is the homology of
$\Gamma \backslash X$, and hence of $\Gamma $.
\item This cell structure is dual to the \Vor \ decomposition in the
following sense:  every $k$-cell of $W$ transversely
intersects exactly one \Vor \ cell of codimension~$k$.
\end{enumerate}
Let $W_{k}$ denote the set of $k$-cells of $W$.  Given $\tau \in W_{k}$, we
denote its dual cell by $\hat \tau \in C^{k}$.  Let $\Z [W_{k}]$
(respectively $\Z [C^{k}]$) denote the free abelian group on the
elements of $W_{k}$ (resp. $C^{k}$).

We may choose orientations compatibly between the cell and cocell
structures in the following sense: for each pair $(\tau ,\hat \tau) $
we may fix orientations so that in the homeomorphism
$$\Int \tau \times \hat \tau  \longrightarrow \R ^{N} $$ 
the product of the orientations is carried to a fixed orientation of
$\R ^{N}$.  This constructs a map $\Z [C^{k}]\rightarrow \Hom _{\Z
}(\Z [W_{k}], \Z )$.

Now to construct $\delta ^{*}$, we use the boundary map from $W_{*}$.
Given two cells $\sigma ,\tau \in W_{*}$, write $\sigma < \tau $ if
$\sigma $ appears in the closure of $\tau $.  Then $\partial
_{k}\colon \Z [W_{k}]\rightarrow \Z [W_{k-1}]$ has the form
\begin{equation*}
\tau \longmapsto \sum _{\sigma < \tau } [\sigma : \tau ]\,\sigma,  
\end{equation*}
where the $[\sigma :\tau ] = \pm 1$ keeps track of the relative
orientation between $\sigma $ and $\tau $.  (Saying $W$ is a regular
cell complex makes $[\sigma :\tau ]$ well defined.)  We define $\delta
^{k}\colon \Z [C^{k}]\rightarrow \Z [C^{k+1}]$ by
\begin{equation}\label{coboundary}
\hat \tau \longmapsto \sum _{\tau < \sigma } [\tau : \sigma  ]\,\hat \sigma 
\end{equation}

\begin{proposition}
With the above coboundary map, $H^{k}(\Gamma ;\Z )$ is naturally
isomorphic to the $k^{{\rm th}}$-cohomology of the quotient modulo
$\Gamma $ of $(\Z [C^{* }],\delta ^{* })$.
\end{proposition}

\begin{proof}
We must show that $\delta ^{2}\equiv 0$ and that $\delta $ is the
adjoint of $\partial $ with respect to the pairing between $\Z
[W_{*}]$ and $\Z [C^{*}]$.  The former is purely formal using
\eqref{coboundary} and the fact that $\hat \tau $ is the map $\Z
[W_{k}]\rightarrow \Z $ that takes $\tau $ to $1$ and all others to
$0$.  The latter is easily verified from the definitions and the
choice of orientations.
\end{proof}

\section{\Vor \ reduction}\label{vor.red} 
In this section we address two questions:
\begin{enumerate}
\item How do we construct $\Pi $ in practice?
\item Given a point $x\in C$, can we determine a
top-dimensional \Vor \ cone containing it?  (Such a cone is unique for generic $x$, and for any given $x$ there
are at most a finite number of such cones containing it.)
\end{enumerate}
In~\cite{voronoi2}, \Vor \ answers these 
in the setting of Example~\ref{classical}, where $C$ is the cone
of real positive-definite symmetric matrices.  In this
section we prove that \Vor 's results remain true in our more general context.

\subsection{}\label{bkpoly}
First we describe some geometric properties of $\Pi $ which are proved
in
\cite{ash}.  

Let $F$ be a facet of $\Pi $, that is, a codimension-one face of $\Pi
$.  Then there is a unique point $y_{F}\in C\cap V(\Q )$ such that 
\begin{enumerate}
\item $F = \left\{\,x\in \Pi \,|\,\langle x,y_{F}\rangle = 1
\right\}$, and 
\item for all $x\in \Pi - F$, we have $\langle x,y_{F}\rangle > 1$.
\end{enumerate}
We say that $y_{F}$ defines a {\em supporting hyperplane} of $\Pi $.  

Given a facet $F$, let $Z_{F}$ be the finite set of points $z\in
L'\cap \Xi (C)$ such that $\langle z,y_{F}\rangle = 1$.  Then $F$ is
the convex hull of $Z_{F}$, and as we range over all $w\in L'\cap \Xi
(C)$ such that $w\not \in Z_{F}$, the set of numbers $\langle
w,y_{F}\rangle$ is bounded below away from $1$.  We call $y_{F}$ the
{\em perfect form} associated to $F$ and $Z_{F}$ the set of {\em
minimal vectors} of $y_{F}$.  In the case of Example \ref{classical},
the $y_{F}$ are perfect quadratic forms in the classical sense, with
minimal vectors $Z_{F}$~\cite{voronoi2}.

Let $\sigma \subset \bar C$ be a rational polyhedral cone.  Then
$\sigma$ satisfies the ``property of Siegel'' with respect to the \Vor
\ decomposition.  Specifically, the intersection $\sigma \cap \Pi $ is
cut out from $\Pi $ by a finite number of supporting hyperplanes
(\cite{ash}, p. 73).  This implies that for any $\sigma $ and for any
$x\in C$, the orbit $\Gamma_{L}x$ meets $\sigma $ in a finite set.

Given any $y\in C (\Q )$, let $\pi (y):V\rightarrow \R $ denote the
linear map $x\mapsto \langle x,y\rangle$.  We also need the following
finiteness result.

\begin{proposition}\label{finiteness.prop}
Let $y\in C (\Q )$.  Then for any $\mu >0$, the set 
\[
\left\{z\in L'\cap \Xi (C)\mid 0< \langle
z,y\rangle\leq \mu \right\}
\]
is finite.
\end{proposition}

\begin{proof}
Given any $\lambda \in \R $, let $H_{\lambda }$ be the affine
hyperplane $\left\{x\in V\mid \langle x,y\rangle = \lambda
\right\}$.  Then $H_{0} = \ker (\pi )$ is rational, since $y\in C (\Q
)\subset V (\Q )$.  Hence the map $\pi (y)$ takes $L$ onto a
lattice in $\R $.  Since some multiple of $y$ lies in $L'$, this
lattice is nontrivial.  Thus to prove the claim, it is enough
to show that for any $\lambda>0$, the set $H_{\lambda }\cap
L'\cap \Xi (C)$ is finite.

To see this, consider the set $\bar C\cap H_{\lambda }$.  Note that
$H_{\lambda }$ meets $C$, since $y' := \displaystyle\frac{\lambda
y}{\langle y,y\rangle}\in H_{\lambda }$.  Let $\ell \subset H_{\lambda
}$ be any line through $y'$, and let $\partial C$ be $\bar C\setminus
C$.  Since $C$ is a cone, $\ell $ must leave $C$, and hence $\ell \cap
\partial C\not = \emptyset $.  Let $x$ be a point in $\ell \cap
\partial C$.  The self-adjointness of $C$ implies that there is
another point $z$ on $\ell \cap \partial C$. (See
Figure~\ref{cap.partial}.)

\begin{figure}[ht]
\begin{center}
\psfrag{0}{$0$}
\psfrag{y}{$y'$}
\psfrag{x}{$x$}
\psfrag{z}{$z$}
\psfrag{l}{$\ell $}
\includegraphics{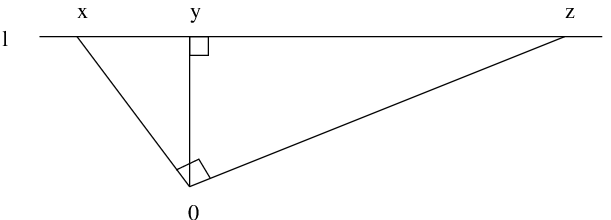}
\caption{\label{cap.partial}}
\end{center}
\end{figure}

Now let $\ell $ range over all lines in $H_{\lambda }$ through $y'$,
and consider the set of lengths of the segments $xy'$ and $y'z$.  This
set is bounded above, and hence $\bar C\cap H_{\lambda }$ is a closed
and bounded subset of $V$.  Since $L$ is a lattice, the set
$H_{\lambda }\cap L'\cap \Xi (C)$, if nonempty, is finite.
\end{proof}
\subsection{}\label{const.pi}
Now we describe the construction of $\Pi $.  Call two facets of $\Pi $
{ \em neighbors} if they meet along a codimension-two face of $\Pi $.
We show how, given a facet $F$, one may systematically find all
neighbors of $F$.

\begin{lemma}\label{neighboring.facets}
Let $F$ and $G$ be neighboring facets of $\Pi $ with perfect forms
$y_{F}$ and $y_{G}$.  Let $v\in V (\Q )$ be orthogonal to the affine
span of the origin and the polytope $F\cap G$, and such that 
$\langle x, v \rangle \geq 0$ for all $x\in F$.
Then $y_{G} = y_{F}+\rho v$ for a unique $\rho \in \R^{>0} $.
\end{lemma}

\begin{proof}
First note that the affine span of $F\cap G$ and the origin is a
hyperplane, since $F\cap G$ is a codimension-two face of $\Pi $.  Thus
$v$ is unique up to a scalar.  Let $v'=y_{G}-y_{F}$.  Then $\langle
x,v' \rangle =0$ for all $x\in F\cap G$, and $v'\not = 0$.  Thus $\rho
v = v'$ for some nonzero $\rho \in \R $, which shows that $y_{G} =
\rho v + y_{F}$.  We must show $\rho $ is positive.

Let $x\in F-(F\cap G)$.  Then $1 < \langle x,y_{G} \rangle = \langle
x,y_{F} \rangle + \rho \langle x,v \rangle = 1 + \rho \langle x, v
\rangle$, which means that $\rho \langle x, v \rangle > 0$.  Since
$\langle x, v \rangle \geq 0$ for all $x\in F$, the result follows.

\end{proof}

Suppose that we are given a facet $F$ with corresponding perfect form
$y_{F}$ and minimal vectors $Z_{F}$.  Choose a maximal proper face $E\subset
F$, and let $Z_{E}\subset Z_{F}$ be the minimal vectors affinely
spanning $E$.  Let $G$ be the facet neighboring $F$ along $E$, and
write $y_{G} = y_{F} + \bar \rho v$, where $v$ is a vector satisfying
the conditions of Lemma~\ref{neighboring.facets}.  Let $Z_{G}$ be the set
of minimal vectors of $y_{G}$.  Define the function $\rho (x)$ by
$$\rho (x) := \frac{1-\langle x, y_{F}\rangle}{\langle x, v\rangle} $$
and 
define $S$ by 
$$S := \left\{x\in L' \cap \Xi (C)\,|\,\langle x, v \rangle < 0 \quad \hbox{and}
\quad y_{F} + \rho (x) v \in C\right\}. $$
Note that $S$ is nonempty.  This follows because whenever $z\in
Z_{G}\setminus Z_{E}$, a computation shows $\rho (z)=\bar \rho $ and
$y_{F} + \rho (z)v = y_{G}$, implying $z\in S$.  This same computation
shows conversely that if $z\in L'\cap \Xi (C)$ and $\rho (z)=\bar \rho $, then
$z$ is a minimal vector of $y_{G}$ not in $Z_{E}$.

\begin{lemma}\label{min.pos.value}
As $x$ ranges over $S$, we have $\rho (x)\geq \bar \rho $.
\end{lemma}

\begin{proof}
Assume there is some
$x'\in S$ with $\rho (x') < \bar \rho $.  Then 
\begin{align*}
\langle x', y_{G} \rangle &= \langle x', y_{F} \rangle + \bar \rho
\langle x', v \rangle\\
	&<\langle x', y_{F} \rangle + \rho (x') \langle x', v \rangle\\
	&=1, 
\end{align*}
which is a contradiction.
\end{proof}

Now we show how to compute $y_{G}$.  Choose any $\ell \in S$ and
consider the point $y_{\ell } := y_{F} + \rho (\ell ) v$.  Let $T$ be
the set
$$T := \left\{x\in L' \cap \Xi (C)\,|\, \langle x, y_{\ell } \rangle \leq 1 \right\}. $$
\begin{lemma}\label{set.tee}
The set $T\cap ( Z_{G}\setminus Z_{E})$ is finite and nonempty.
\end{lemma}

\begin{proof}
Since $v\in V (\Q )$, we have that $y_{\ell }\in C (\Q )$.  Thus $T$
is finite by
Proposition~\ref{finiteness.prop}.

We now show $T\cap ( Z_{G}\setminus Z_{E})$ is nonempty.  Let $z\in (
Z_{G}\setminus Z_{E})\cap S$, a set we have seen is nonempty.  Then
\begin{eqnarray*}
\langle z, y_{\ell } \rangle &=& \langle z, y_{F} \rangle + \rho(\ell )\langle z, v \rangle\\
&\leq& \langle z, y_{F} \rangle + \bar \rho\langle z, v \rangle\quad
\hbox{[Lemma~\ref{min.pos.value}]}\\
&=& \langle z, y_{G} \rangle \\
&=& 1.
\end{eqnarray*}
Thus $z\in T$, and $T\cap ( Z_{G}\setminus Z_{E})$ is nonempty.
\end{proof}

\begin{theorem}\label{finite.num.of.steps}
Given a facet $F$ of the \Vor \ polyhedron, all its neighbors may be
determined in a finite number of steps.
\end{theorem}

\begin{proof}
Let $y_{F}$ be the perfect form corresponding to $F$, and let $Z_{F}$
be the set of minimal vectors of $y_{F}$.  Using standard techniques
of convex geometry, such as Fourier-Motzkin elimination
(\cite{ziegler}, p. 37), we may determine all the maximal proper faces of
$F$.  Let $E$ be such a face, and let $Z_{E}\subset Z_{F}$ be the
minimal vectors affinely spanning $E$.  Using $Z_{E}$, we may
determine $v$ as in Lemma~\ref{neighboring.facets}.  Let $G$ be the
facet of $\Pi $ neighboring $F$ along $E$, and let $y_{G}$ be the
corresponding perfect form.  We need to compute $y_{G}$.

First we find an $\ell \in S$ by searching over $L$, and then using
$\ell $ we construct the finite set $T$.  The latter can done in a
finite number of steps because it is equivalent to finding the set of
vectors in a lattice on which a positive-definite quadratic form is
less than a constant (cf. Example~\ref{classical} and the final
paragraph of \S\ref{bkbegin}).  By Lemma~\ref{set.tee}, $T\cap (
Z_{G}\setminus Z_{E})$ is nonempty, where $Z_{G}$ is the set of
minimal vectors of $y_{G}$.  Now let $Z\subset ( T\setminus Z_{E})$ be
the set on which $\rho (x)$ attains its minimum.  By
Lemma~\ref{min.pos.value} and the paragraph preceding it, $Z\subset
Z_{G}\setminus Z_{E}$.  Let $H$ be the affine span of $Z\cup Z_{E}$.
Then $y_{G}$ is the unique point satisfying $\langle x,y_{G}\rangle =
1$ for all $x\in H$.

Repeating this procedure for each maximal face of $F$, we may
determine all the neighbors of $F$.
\end{proof}

\begin{remark}\label{about.s}
In practice, it may be the case that $S$ \emph{only} consists of
$Z_{G}\setminus Z_{E}$, as a computation with $SL_{2} (\Z )$ shows.    
\end{remark}

Hence one may find facets of $\Pi $ provided one can construct an
initial facet.  In the setting of Example~\ref{classical}, \Vor \ did
this by showing that the quadratic form $A_{n}$, defined by
$$\sum _{i=1}^{n} x_{i}^{2} + \sum _{1\leq i<j\leq n}
(x_{i}-x_{j})^{2},$$ 
is perfect for all $n$ (\cite{voronoi2}, \S 29).\footnote{\Vor \
called 
this form the \emph{principal perfect form}.}  In our
more general setting, one cannot write down a perfect form that works
for every case, even if one restricts to Bianchi groups (Example
\ref{imag.quad}).  However, in practice one may do the following.

First choose a large bounded set $U\subset V$ containing the
origin, and let $\Sigma $ be the convex hull of $L'\cap \,\Xi (C)\cap
U$.  Then $\Sigma$ is a bounded polytope in $V$.  Furthermore, since
the facets of $\Pi $ are bounded, if $U$ is sufficiently large
then most facets of $\Sigma $ will be facets of $\Pi $.  

To check if a facet $F\subset \Sigma $ is a facet of $\Pi $, one
computes $y_{F}$ and checks whether the $z\in L'\cap \,\Xi (C)$ such
that $\langle z,y_{F} \rangle = 1$ are vertices of $F$.  

\begin{remark}
In the Bianchi case, another possibility is to first construct the
retract $W$ from \S\ref{connection} using techniques
in~\cite{vogtmann}, and then use the duality between them to deduce
the structure of $\Pi $.
\end{remark}

\subsection{}\label{vor.alg}
To answer question (2), \Vor \ describes a reduction algorithm.  This
algorithm is based on the following:

\begin{proposition}\label{finite.beta}
Fix an $x\in \Pi $ and a real number $\mu \geq 1$.  Then there are only a
finite number of perfect forms $y_{F}$ satisfying 
$$ \langle x,y_{F}\rangle \leq \mu . $$
\end{proposition}

\begin{proof}
Choose a facet $F$ such that $\langle x,y_{F}\rangle \leq \mu$.  Since
$y_{F}\in C (\Q )$, Proposition~\ref{finiteness.prop} implies the set
$\left\{ z\in L'\cap \Xi (C) \, | \, \langle z,y_{F} \rangle \leq \mu
\right\}$ is finite.  Hence the polyhedron
$$\Sigma = \left\{\,y\in \Pi \,|\, \langle y,y_{F}\rangle \leq
\mu\,\right\} $$ 
is a bounded polytope in $\bar C$.  If $\mu \in \pi (y_{F}) (L)$,
where $\pi $ is defined in \S\ref{bkpoly}, then the
hyperplanes bounding $\Sigma $ will be rational with respect to the
$\Q $-structure on $V$, and thus $\Sigma $ will have vertices in $V(\Q
)$.  Hence, replacing $\mu $ by a slightly larger number if
necessary, we may assume $\Sigma $ has vertices in $V(\Q )$.

Therefore the cone generated by the vertices of $\Sigma $ is rational
polyhedral, and only meets a finite subset of the orbit $\Gamma_{L}x$.
By taking the adjoint action of $\Gamma_{L}$ with respect to
the scalar product, it follows that only finitely many facets $F'$
that are $\Gamma_{L}$-equivalent to $F$ satisfy $ \langle
x,y_{F'}\rangle \leq \mu$.  Since there are only a finite number of
facets modulo $\Gamma_{L}$, the result follows.  
\end{proof}

The following lemma gives a local condition for when a given point of $C $
lies in a cone over a facet of $\Pi $.

\begin{lemma}\label{local.cond}
Let $F$ be a facet of $\Pi $, and let $\G$ be the set of neighbors of
$F$.  Let $x\in C$.  Suppose that $\langle
x,y_{F}\rangle \leq \langle x,y_{G}\rangle $ for all $G\in \G$. 
Then $x$ lies in the cone over $F$.
\end{lemma}

\begin{proof}
Choose $\lambda $ so that $x' := \lambda x$ satisfies $\langle
x',y_{F}\rangle = 1$ and $\langle
x',y_{G}\rangle \geq 1$ for $G\in \G$.  For $\epsilon >0$, let $\Sigma _{\epsilon }$ be
the polyhedron defined by 
$$\Sigma _{\epsilon } = \left\{ x \,|\, 1\leq \langle x, y_{F} \rangle
\leq 1+\epsilon \quad \hbox{and} \quad \langle x, y_{G} \rangle \geq 1 \quad \hbox{for $G\in \G$.}\right \} $$
If $\epsilon $ is sufficiently small, then $\Sigma _{\epsilon }\subset
\Pi $.  Hence $x'\in \Sigma _{\epsilon }$, and thus $x'\in \Pi $.
Since $x'$ also lies in the supporting affine hyperplane $\{x \,| \,\langle
x,y_{F} \rangle = 1\}$, it must lie in $F$.
\end{proof}

Now we describe our \Vor \ reduction algorithm.

\begin{theorem}\label{reduction.alg}
Let $x\in \Pi $, and choose a facet $F$.  Let $\mu =\langle
x,y_{F}\rangle$.  The following algorithm determines a cone in the
\Vor \ decomposition containing $x$:
\begin{enumerate}
\item For each neighbor $G$ of $F$, compute $\langle
x,y_{G}\rangle $.  
\item If there exists a neighbor $G$ with $\langle
x,y_{G}\rangle < \mu $, replace $F$ with $G$, $\mu $ with $\langle x,
y_{G} \rangle$, and return to step one.
\item Otherwise, terminate the procedure:  $x$ lies in the cone
generated by $F$.
\end{enumerate}
\end{theorem} 

\begin{proof}
By Lemma~\ref{local.cond}, if the algorithm terminates then we have
determined a cone containing $x$.  We now prove that the algorithm
terminates.  Suppose that a neighbor $G$ of $F$ has $\langle
x,y_{G}\rangle < \mu $. (Note that this quantity must be positive
since $C$ is self-adjoint and $x,y\in C$.)  Then we return to step
one, and we have decreased the scalar product.  Since by Proposition
\ref{finite.beta} the set of facets satisfying $\langle
x,y_{F'}\rangle \leq \mu $ is finite, the algorithm must terminate.
\end{proof}

\begin{remark}\label{comp.complexity}
The data needed to implement this algorithm is the same as that needed
for the structure of $\Pi $ modulo $\Gamma _{L}$, along with some additional
information.  In particular, one must determine: 
\begin{enumerate}
\item A finite set $\mathscr{F}$ of representatives of the facets of $\Pi $ modulo
$\Gamma _{L}$. 
\item For each $F\in \mathscr{F}$ and each neighbor $G$ of $F$, a
element $\gamma \in \Gamma _{L}$ such that $\gamma G \in \mathscr{F}$.
\end{enumerate}
For an example of the
implementation of this algorithm for $SL_{2}(\Z )$, we refer to
\cite{jaquet}.

As far as we know, the computational complexity of this
algorithm is unknown.  However, in our experience it performs
very well for $SL_{3}(\Z )$ and $SL_{4}(\Z )$.
\end{remark}

\section{The modular symbol algorithm}\label{mod.sym}
In this section we define modular symbols and describe our Hecke
algorithm.  Our definition of a modular symbol is closely related to the
definitions appearing in~\cite{ash.minmod} and~\cite{ash.rudolph}.  As
before, let $N$ be the dimension of $X$ and let $d$ be the 
cohomological dimension of $\Gamma $.  We assume from now on that $G$
has $\Q $-rank one, so that $N=d+1$.

\subsection{}\label{mod.symbols}
Let $\rho (\Pi )$ be the set of rays in $\bar C$ generated by the
vertices of the \Vor \ polyhedron $\Pi $.  Given an ordered pair
$(u,v)\in \rho (\Pi )\times \rho (\Pi )$, we want to construct a class
$[u,v]\in H^{d}(\Gamma ;\Z )$. 

To this end we recall that $X$ may be extended to a bordification
$\bar X$ such that the quotient $\Gamma \backslash \bar X$, the
\emph{Borel-Serre compactification}, is a compact manifold with
corners with interior $\Gamma \backslash X$~\cite{borel.serre}.  Let
$\pi \colon \bar X \rightarrow \Gamma \backslash \bar X$ be the
canonical projection.  Given $u,v\in \rho (\Pi) $, we determine a path
in $\Gamma \backslash \bar X$ as follows.  First $u$ and $v$ determine
a closed cone $\sigma \subset V$, and we let $\bar \sigma$ be the
closure of $\sigma \cap V$ mod homotheties in $\bar X$.  Then $\pi
(\bar \sigma )$ is a path with endpoints lying in $\partial (\Gamma
\backslash \bar X)$.  Choosing an ordering $(u,v)$ fixes an
orientation of $\pi (\bar \sigma )$ and thus determines a class in
$H_{1}(\Gamma \backslash \bar X, \partial (\Gamma \backslash \bar X);
\Z )$.  By Lefschetz duality,
$$ H_{1}(\Gamma \backslash \bar X, \partial (\Gamma \backslash \bar X);
\Z ) = H^{d}(\Gamma \backslash \bar X;\Z ),$$
and since $\Gamma \backslash \bar X$ is homotopy equivalent to $\Gamma
\backslash X$, we have actually determined a class in $H^{d}(\Gamma
;\Z )$. 

\begin{definition}
A {\em modular symbol} is a class in $H^{d}(\Gamma ;\Z )$ 
constructed as above from an ordered
pair $(u,v)\in \rho (\Pi )\times \rho (\Pi )$.  The class is denoted $[u,v]$.
\end{definition}
 
Note that this definition is almost the same as that in \S\ref{manin.work} for
$\Gamma \subset SL_{2}(\Z )$, because using $(u,v)$ to generate a
cone is the same as choosing a specific path as in \S\ref{manin.work}.
However, as in the earlier definition, the class $[u,v]$ is
independent of this path.  

We have the following analogue of Proposition~\ref{manin.sl2}:

\begin{proposition}\label{properties}
Let $u,v\in \rho (\Pi )$.  The modular symbols satisfy the following:
\begin{enumerate}
\item $[u,v]=-[v,u]$.
\item If $w\in \rho (\Pi )$, then $[u,v]=[u,w]+[w,v]$.
\item The modular symbols span $H^{d}(\Gamma ;\Z )$.
\end{enumerate}
\end{proposition}

\begin{proof}
Only (2) and (3) require proof.  To prove (2), let $\sigma $ be the
cone generated by $u$ and $v$, and choose a ray $x\subset \sigma $
distinct from $u$ and $v$.  Let $\phi \colon [0,1]\rightarrow \bar
C$ be a continuous family of rays such that $\phi (0) = x$ and $\phi
(1)=w$.  Let $\sigma _{1} (t)$ (respectively $\sigma _{2} (t)$) be the cone
generated by $u$ and $\phi (t)$ (resp. $\phi (t)$ and $v$).  Then
$\phi $ provides a continuous deformation of $\sigma $ into $\sigma
_{1} (1)\cup \sigma _{2} (1)$ that induces the homology   $[u,v]=[u,w]+[w,v]$.

Now to prove (3), note that the results of \S\ref{connection} imply
that any class in $H^{d}(\Gamma ;\Z )$ may be written as a cocycle for
$\delta$ using the \Vor \ cones of codimension-$d$ modulo $\Gamma $.
But the \Vor \ cones of codimension $d$ are cones generated by pairs
of vertices of $\Pi $, and so any such cycle is in the span of the
modular symbols.
\end{proof}

\begin{remark}
The third statement is also a consequence of the results in~\cite{ash.minmod}.
\end{remark}

\subsection{}
Now we describe our replacement for the notion of a unimodular symbol.
Recall that $R(\sigma )$ is the set of spanning rays for a cone
$\sigma $ (\S\ref{bkend}).

\begin{definition} Let $\sigma \subset \bar C$ be a rational
polyhedral cone.  Then $\sigma $ is {\em \Vor -reduced} if there is a
top-dimensional \Vor \ cone $\sigma _{\alpha }$ such that
$$R(\sigma )\subset R(\sigma _{\alpha }). $$
\end{definition}

Note that a \Vor -reduced cone is spanned by cusps.
A \Vor -reduced cone need not be a \Vor \  cone, since the
top-dimensional \Vor \ cones are not simplicial in general (cf. Examples
\ref{classical} and~\ref{imag.quad}).  However, since every \Vor \
cone is generated by a finite set of cusps, and because of 
the finiteness properties of any $\Gamma $-admissible decomposition,
we have 

\begin{proposition}\label{modgamma}
Modulo $\Gamma $, there are only finitely many \Vor -reduced cones.
\end{proposition}

We say that a modular symbol is \Vor -reduced if it is induced by a
\Vor -reduced cone.  Propositions~\ref{properties} and 
\ref{modgamma} imply that the \Vor -reduced modular symbols
provide a finite spanning set for $H^{d}(\Gamma ;\Z )$.

\begin{remark}\label{spanning.set}
Although the \Vor -reduced modular symbols provide a finite spanning
set for $H^{d}(\Gamma ;\Z )$, they are not a basis of $H^{d}(\Gamma
;\Z )$.  In fact, they are not even a basis of $\Z [C^{*}]$
modulo $\Gamma $, because they are not necessarily supported on
codimension-$d$ \Vor \ cones.  However, in practice this does not
affect their usefulness (cf. Remark~\ref{precompute}).

\end{remark}

\subsection{}\label{algorithm}
Let $u,v\in \rho (\Pi )$, and let $\sigma $ be the rational
polyhedral cone generated by $u$ and $v$.  We are now ready to
describe and prove our algorithm.

Given $x\in C$, let $\sigma (x)$ be the unique \Vor \ cone containing
$x$, and by abuse of notation let $R(x)$ be $R(\sigma (x))$. 

Let $x_{1},\ldots,x_{r}$ be points in $\sigma \cap C$ such that the
rays $\R^{>0}x_{1},\ldots,\R^{>0}x_{r}$ are distinct.  These points
subdivide $\sigma $ into a collection of cones, namely those generated
by the pairs $(u,x_{1}), \ldots,(x_{r},v)$.  

\begin{definition}\label{suff.fine}
The decomposition of $\sigma $ by the $x_{i}$ as above is called a {\em
sufficiently fine partition} if 
\begin{enumerate}
\item $u\in R(x_{1})$,
\item $v\in R(x_{r})$, and 
\item for $i=1,\ldots,r-1$, we have $R(x_{i})\cap R(x_{i+1})\not
=\emptyset $. 
\end{enumerate}
\end{definition}

\begin{lemma}
Sufficiently fine partitions of $\gamma $ exist.
\end{lemma}

\begin{proof}
Since $\sigma $ is rational polyhedral, the Siegel property implies
that $\sigma \cap \Pi $ is cut out by finitely many supporting
hyperplanes.  Hence $\sigma \cap C$ meets only finitely many
top-dimensional cones, and the intersection of these cones with
$\sigma $ subdivides the latter into finitely many $2$-cones $\{\sigma
_{i} \}$.  We may take $x_{i}$ to be any nonzero point in the interior
of $\sigma _{i}$.  Conditions (1) and (2) of Definition
\ref{suff.fine} are trivially satisfied.  Also, $R (x_{i})\cap R
(x_{i+1})\not =\emptyset $ because $\sigma _{i}\subset V_{i}$ and
$\sigma _{i+1}\subset V_{i+1}$, where $V_{i}$ and $V_{i+1}$ are \Vor \
cones that have a face in common.

\end{proof}

\begin{remark}
For computational purposes, one may construct sufficiently fine
partitions of $\gamma $ as follows.  Let $u,v\in \rho (\Pi )$ and let
$\sigma $ be the cone generated by $u$ and $v$.  Choose points $\bar
u\in u$ and $\bar v\in v$.  Let $\bar x$ be the midpoint of the
segment between $\bar u$ and $\bar v$, and let $x$ be the cone
generated by $\bar x$.  Now apply Theorem~\ref{reduction.alg} to check
whether $u\in R (\bar x)$ and $v\in R (\bar x)$.  If these
conditions are not satisfied, bisect the segments between $\bar u$,
$\bar x$ and $\bar x$, $\bar v$, and check conditions (1), (2), and
(3).  Eventually, by the Siegel property, after a finite number of
iterations one will have constructed a sufficiently fine partition of
$\sigma $.
\end{remark} 

Now we present our algorithm.

\begin{theorem}\label{hecke.alg}
Given a modular symbol $[u,v]$, the following constructs a chain of 
\Vor -reduced modular symbols homologous to $[u,v]$:
\begin{enumerate}
\item Choose a set of points $\{ x_{i}\}$ inducing a sufficiently fine
partition of the cone generated by $u$ and $v$.  
\item For $i=1,\ldots,r-1$, choose a ray $q_{i}\in R(x_{i})\cap
R(x_{i+1})$.   
\end{enumerate}
Then $[u,v] = [u,q_{1}] + [q_{1},q_{2}] + \cdots + [q_{r},v]$.
\end{theorem} 

\begin{proof}
First note that each modular symbol on the right hand side is \Vor
-reduced, since $q_{i}$ and $q_{i+1}$ are both rays from
$R(x_{i+1})$.   We must
show there is a homology between the right side and the left.
Notice that 
$$[u,v] = [u,q_{1}] + [q_{1},v] $$
by Proposition~\ref{properties}.  Repeatedly applying this proposition, we see that 
$$[q_{i},v] = [q_{i},q_{i+1}] + [q_{i+1},v], $$
which completes the proof.
\end{proof}
 
\begin{remark}\label{precompute}
For computational purposes, to determine the action of a Hecke
operator we must write any modular symbol in terms of a basis of
$H^{d}(\Gamma ;\Z )$, and by Remark~\ref{spanning.set} the technique
in Theorem~\ref{hecke.alg} is not sufficient to do this.  However, in
practice we may precompute explicit homologies between \Vor -reduced
modular symbols and modular symbols supported on \Vor \ cones, as
follows.

Let $\mathscr{F}$ be a set of representatives of the facets of $\Pi $
modulo $\Gamma _{L}$.  For each $F\in \mathscr{F}$, let $u,v$ be any two
vertices of $F$.  Then $[u,v]$ is a \Vor -reduced modular symbol.  To
write $[u,v]$ in terms of the basis of codimension-$d$ \Vor \ cones,
choose any sequence of vertices $u=u_{0},u_{1},\ldots,u_{k}=v$ of $F$
such that $u_{i}$ and $u_{i+1}$ are joined by an edge of $F$.  Then
$$[u,v]=[u,u_{1}]+\cdots +[u_{k-1},v] $$
is the desired homology.  Now repeat for all $u,v$ and all $F\in \mathscr{F}$.
\end{remark}
%
%
\bibliographystyle{amsplain}
\bibliography{msa}

\end{document}